\documentclass[12pt]{article}
\usepackage{ifpdf}
\usepackage {graphics}
\usepackage{amsmath,mathptmx,amssymb,bm}
\newtheorem{theorem}{Theorem}

\begin{document}


\title{Symmetry reductions of a generalized Kuramoto-Sivashinsky equation via equivalence transformations}


\author{R. de la Rosa${}^{a}$, M.S. Bruz\'on${}^{b}$\\
 ${}^{a}$ Universidad de C\'adiz, Spain  (e-mail: rafael.delarosa@uca.es). \\
 ${}^{b}$ Universidad de C\'adiz, Spain  (e-mail:  m.bruzon@uca.es). \\  
}

\date{}
 
\maketitle

\begin{abstract}

In this paper we consider a generalized Kuramoto-Sivashinsky equation. The equivalence group of the class under consideration has been constructed. This group allows us to perform a comprehensive study and a clear and concise formulation of the results. We have constructed the optimal system of subalgebras of the projections of the equivalence algebra on the space formed by the dependent variable and the arbitrary functions. By using this optimal system, all nonequivalent equations admitting an extension by one of the principal Lie algebra of the class under consideration can be determined. Taking into account the additional symmetries obtained we reduce some partial differential equations belonging to the class into ordinary differential equations. We derive some exact solutions of these equations.\\

\end{abstract}



\maketitle

\section{Introduction}
\label{}

\noindent The Kuramoto-Sivashinsky (KS) equation
\begin{equation}\label{ks1}u_t+\beta_1 u_{xx}+\beta_2 u_{xxxx}+\delta u_{xxx}+uu_x=0,\end{equation} arises in the modelling of flow on an inclined plane, where $\beta_1$,  $\beta_2$ and $\delta$ are constants.  This equation was introduced by Kuramoto \cite{KS1} in
one-spatial dimension for the study of phase turbulance in the Belousov-Zhabotinsky reaction. Sivashinsky derived it independently in the context of small thermal diffusive instabilities for laminar flame fronts. Some authors have previously studied this equation \cite{NK1,NK2}.\\

The study of dynamical behaviours such as chaos arises as an effective way to explain the physics of low-dimensional dynamical systems, i.e., dynamical systems with small number of state variables. Nevertheless most problems and phenomena in physics and other areas of science are described by partial differential equations. Among them we highlight the KS equation. Khellat and Vasegh \cite{khellat} improved the understanding of the connection between low-dimensional systems and the KS equation which is subject to spatially periodic boundary condition and with arbitrary periodic initial condition given by

\begin{equation}\label{ksperiodic}
\begin{array}{l}
u_t+ u_{xx}+k u_{xxxx}+uu_x=0, \qquad (x,t) \in \mathbb{R} \times \mathbb{R}^{+}, \\

u(x,0)=f(x), \qquad u(x-\pi,t)=u(x+\pi,t), \qquad \displaystyle \int_{-\pi}^{\pi} f(x) dx=0,
\end{array}
\end{equation}

\noindent where $k$ represents the viscosity. The authors showed that the KS equation (\ref{ksperiodic}) leads respectively to stability, pitchfork bifurcation and a new type of behaviour. Sahoo and Saha Ray \cite{sahoo} considered the time-fractional KS equation given by
\begin{equation}\label{ksfractional}
D_t^{\alpha} u+ b u_{xx}+ k u_{xxxx}+a uu_x=0,\\
\end{equation}

\noindent which is a generalization of equation (\ref{ksperiodic}), where $0<\alpha \leq 1$ and $a$, $b$ and $k$ are arbitrary constants. The authors obtained new types of exact analytical solutions of equation (\ref{ksfractional}) by applying the tanh-sech method with the help of the fractional complex transform. \\

\noindent There have been several generalizations
of the KS equation such as the generalized KS with
dispersive effects \cite{guo}

\begin{equation}\label{edg} u_t+\left( f(u) \right)_x+\alpha
u_{xx}+\phi(u)_{xx}+\beta u_{xxx}+\gamma
u_{xxxx}=g(u),\end{equation} where $\alpha, \beta$ and $\gamma$
are constants, $f(u),$ $\phi(u)$ and $g(u)$ are functions of the dependent variable $u$. Here, $\left( f(u) \right)_x= \frac{\partial f(u)}{\partial u} u_x= f_u u_x$.
In \cite{BGC} Bruz\'on {\it et al.} considered equation (\ref{edg}) and they made a full
analysis of the symmetry reductions and  proved that the nonclassical method applied to the equation leads to new reductions which cannot be obtained by Lie classical symmetries. Furthermore, making use of the new similarity solutions that the nonclassical method determines, they obtained some exact solutions which cannot be determined by applying a Lie classical reduction. If $\alpha$ and $\beta$ are constants, $\gamma=1$, $f(u)= \frac{\alpha_1}{2}u^2+\beta u+\alpha_5$, $g(u)=0$ and
$\phi(u)=\frac{\alpha_2}{2}u^2+(\alpha_3-\alpha)u+\alpha_4$, with $\alpha_i$, $i=1,\ldots,5$, arbitrary constants, equation (\ref{edg}) includes the Korteweg-de Vries equation
supplements by additional terms of the KS equation which describes nonlinear convection and the input of energy produced by Marangoni forces on
the long scales together with energy dissipation on short scales. For this equation Bruz\'on {\it et al.} obtained classical and nonclassical symmetries and the reduced equations \cite{app}. We propose herein a generalized KS equation
\begin{equation}
\label{ed1} u_t+ \left(f(u)\right)_x+\alpha(u) u_{xx}+\left(\phi(u)\right)_{xx}+\beta(u) u_{xxx}+\gamma(u) u_{xxxx}= g(u),
\end{equation}

\noindent with $f(u)$, $g(u)$, $\alpha(u)$, $\beta(u)$, $\gamma(u) \neq 0$ and $\phi(u)$ arbitrary functions. Here $\left( \phi(u) \right)_{xx}= \frac{\partial^2 \phi(u)}{\partial u^2} u_x^2+ \frac{\partial \phi(u)}{\partial u} u_{xx}= \phi_{uu} u_x^2+\phi_u u_{xx}$.\\

\noindent Since in different fields arise many equations that involve several arbitrary functions, such as
industrial mathematics, mathematical biology and fluid mechanics, the study of variable-coefficient equations has increased in the last decades \cite{Bru:14,bgi,RGB:15,FRE:14,GG,JoKh:10,JoKh:11, MoRa:12, Tra:14, Tra:16}.\\

 On the other hand, equivalence transformations are used mainly for classifying classes of differential equations with arbitrary functions \cite{RGB:16a,RGB:16b,GI,SM,ovsian,RP,senthil,sophorita2008}. Just like symmetries of a differential equation transform solutions of the equation into other solutions, equivalence transformations map every equation of a class $\mathfrak{C}$ of differential equations into an equation of the same class. Ovsiannikov \cite{ovsian} defined a methodology and notation for dealing with such transformations, for which he used the term equivalence transformations. He derived some results about them, including important properties:
\begin{enumerate}
\item[a)] The transformations act on every equation in the class $\mathfrak{C}$.
\item[b)] The transformations are fixed point transformations, in the sense that they
do not depend on the arbitrary elements, and are realized on the point
space (independent and dependent variables) associated with the differential
equations.
\item[c)] The transformations act on the arbitrary elements as point transformations
of an augmented space of independent and dependent variables and additional
variables representing values taken by the arbitrary elements.
\end{enumerate}
The collection of all such transformations constitutes Ovsiannikov's equivalence
group. In \cite{AGI} Akhatov, Gazizov and Ibragimov, by using Ovsiannikov's method, determined the infinitesimal
form of transformations for a potential form of the nonlinear diffusion equation.\\

 Ibragimov, Torrisi and Valenti used the equivalence group to give a preliminary symmetry group classification. They found the equivalence group for a large class of nonlinear hyperbolic equations and executed the preliminary classification for a finite-parameter subgroup \cite{ITV}. By using the weak equivalence classification, Torrisi and Tracin\`a \cite{TT:98} obtained a classification with respect to thermal
conductivity and the internal energy at equilibrium of a system of partial differential equations (PDEs) which describe the unidimensional heat conduction in a homogeneous isotropic rigid body. In \cite{RT}, Lie symmetries for a class of drift-diffusion systems were found following a different procedure based on the weak equivalence classification already introduced
in \cite{TT:98,TTV:94,TTV:96}. They authors claimed that {\it even if this approach does not guarantee obtaining a complete symmetry classification, yet this method provides a systematic way
for obtaining wide classes of symmetries when arbitrary functions appear}.\\

 In \cite{VPS:14}, Vaneeva, Popovych and Sophocleous discussed how point transformations can be used for the study of integrability, in particular,
for deriving classes of integrable variable-coefficient differential equations. They described the procedure of
finding the equivalence groupoid of a class of differential equations focusing on the case of evolution equations and they applied this study to a class of fifth-order variable-coefficient
Korteweg-de Vries-like equations. Two alternative ways to solve completely the group classification problem for a variable coefficient Gardner equation were presented in \cite{Sopho:14}: the gauging of arbitrary elements of the class by the equivalence transformations and the method of mapping between classes. Recently, Tracin\`a \cite{Tra:15} proved a criterion, based on the property of nonlinear self-adjointness, for the existence of an invertible point transformation which maps a given PDE to a linear PDE.\\

On the other hand, Gandarias and Ibragimov \cite{GI} found and employed the Lie algebra of the generators of the equivalence transformations of a fourth-order non-linear evolutionary equations
given by
$$\begin{array}{c}
 u_t + f(u) u_{xxxx} + g(u) u_x u_{xxx} + h(u) u_{xx}^2 +
 d(u) u_x^2 u_{xx} - p(u) u_{xx} - q(u) u_x^2=0.
\end{array}$$

In this paper we obtain the equivalence transformations of the generalized KS equation (\ref{ed1}). We determine the generators of the equivalence algebra. Equivalence algebra is used to perform a preliminary group classification. By using a theorem on projections we solve the problem of preliminary group classification of class (\ref{ed1}) by means of the construction of the optimal system of subalgebras of the nonzero projections of the equivalence algebra on the space $\left( u,f,g,\alpha,\beta,\gamma, \phi \right)$. The optimal system allows us to obtain all nonequivalent equations (\ref{ed1}) admitting an extension by one of the principal Lie algebra of class (\ref{ed1}). Taking into account some nonequivalent equations which admit an additional symmetry, we have reduced equation (\ref{ed1}) with the corresponding coefficients $f,g,\alpha,\beta,\gamma,\phi$ into ODEs. Finally, some exact solutions are obtained.\\

\section{Equivalence transformations}\label{equivalence transformations}

In this section we determine the equivalence transformations of class (\ref{ed1}). An equivalence transformation of class (\ref{ed1}) is a nondegenerate point transformation from $(t,x,u)$ to $(\tilde{t},\tilde{x},\tilde{u})$ with the property that it preserves the differential structure of the equation, that is, it transforms any equation of class (\ref{ed1}) to an equation from the same class but with different arbitrary elements, $\tilde{f}(\tilde{u})$, $\tilde{g}(\tilde{u})$, $\tilde{\alpha}(\tilde{u})$, $\tilde{\beta}(\tilde{u})$, $\tilde{\gamma}(\tilde{u})$ and $\tilde{\phi}(\tilde{u})$ from the original ones. The set of equivalence transformations forms a group denoted by ${\cal E}$.\\

The equivalence transformations for our equation can be obtained by applying Lie's infinitesimal criterion \cite{ovsian}. Unlike classical symmetries, we require not only the invariance of class (\ref{ed1}) but also the invariance of the auxiliary system
\begin{equation}\label{aux}
f_t = f_x=g_t=g_x=\alpha_t=\alpha_x=\beta_t=\beta_x=\gamma_t=\gamma_x=\phi_t=\phi_x=0.
\end{equation}

\noindent We consider the one-parameter group of equivalence transformations in the augmented space $\left( t,x,u,f,g,\alpha,\beta,\gamma, \phi \right)$ given by
\begin{equation}\label{trans}\begin{array}{rcl} \tilde{t} & = & t+\epsilon \, \tau(t,x,u)+O(\epsilon^2),
\\ \tilde{x} & = & x+\epsilon \, \xi(t,x,u)+O(\epsilon ^2),\\ \tilde{u} & = & u+\epsilon \,
\eta(t,x,u)+O(\epsilon ^2),

\\ \tilde{f} & = & f+\epsilon \,
\omega^1( t,x,u,f,g,\alpha,\beta,\gamma, \phi )+O(\epsilon ^2),
\\ \tilde{g} & = & g+\epsilon \,
\omega^2( t,x,u,f,g,\alpha,\beta,\gamma, \phi )+O(\epsilon ^2),
\\ \tilde{\alpha} & = & \alpha+\epsilon \,
\omega^3( t,x,u,f,g,\alpha,\beta,\gamma, \phi )+O(\epsilon ^2),
\\ \tilde{\beta} & = & \beta+\epsilon \,
\omega^4( t,x,u,f,g,\alpha,\beta,\gamma, \phi )+O(\epsilon ^2),
\\ \tilde{\gamma} & = & \gamma+\epsilon \,
\omega^5( t,x,u,f,g,\alpha,\beta,\gamma, \phi )+O(\epsilon ^2),
\\ \tilde{\phi} & = & \phi+\epsilon \,
\omega^6( t,x,u,f,g,\alpha,\beta,\gamma, \phi )+O(\epsilon ^2),
\end{array}
\end{equation} where  $\epsilon$ is the
group parameter. In this case, the vector field takes the following form\\
\begin{equation}\label{generator}
 Y  =   \tau \partial_t + \xi \partial_x+\eta \partial_u+ \omega^1 \partial_f  + \omega^2 \partial_g + \omega^3 \partial_{\alpha} + \omega^4 \partial_{\beta} + \omega^5 \partial_{\gamma}+ \omega^6 \partial_{\phi}.   \end{equation}

\noindent We notice that the fourth derivative with respect to $x$ appears in the equation, therefore we need to take into account the fourth prolongation of the operator (\ref{generator})

\begin{equation}\label{prolong}
 \tilde{Y}  =   Y+ \zeta^t \partial_{u_t}+ \zeta^x \partial_{u_x} +\zeta^{xx} \partial_{u_{xx}} +\zeta^{xxx} \partial_{u_{xxx}}+\zeta^{xxxx} \partial_{u_{xxxx}}+ \tilde{\omega}^i_t \partial_{f^i_t}+ \tilde{\omega}^i_x \partial_{f^i_x}+ \tilde{\omega}^i_u \partial_{f^i_u},   \end{equation}

\noindent where $f^i$, $i = 1,\ldots, 6$, represents each component $\left( f,g,\alpha,\beta,\gamma, \phi \right)$. The  coefficients $\zeta^J$ is defined by

 $$\zeta^J(t,x,u^{(4)})=D_J(\eta-\tau u_t-\xi u_x)+\tau u_{Jt} +\xi u_{Jx},$$

\noindent with $J=(j_1,\ldots,j_k)$,  $1\leq j_k\leq 2$ and $1\leq k\leq 4$, $u^{(4)}$ denotes the sets of partial derivatives up to fourth order and $D_t$, $D_x$ are the total derivatives with respect to $t$ and $x$ \cite{olver}. Lastly, the coefficients $\tilde{\omega}^i_r$ are given by

$$ \tilde{\omega}^i_r= \tilde{D}_r ( \omega^i ) - f_t^i \, \tilde{D}_r \left( \tau \right) - f_x^i \, \tilde{D}_r \left( \xi \right)- f_u^i \, \tilde{D}_r \left( \eta \right),$$

\noindent where

$$ \tilde{D}_r = \partial_r + f^i_r \partial_{f^i}.$$

\noindent For further information on how prolongations
of higher order can be obtained, one can refer to references \cite{senthil,sophorita2008}.\\

\noindent The invariance of system (\ref{ed1})-(\ref{aux}) under the one-parameter group of equivalence transformations (\ref{trans}) with infinitesimal generator (\ref{generator}) leads to a system of determining equations.  After having solved the determining system, omitting tedious calculations, we obtain the associated equivalence algebra $\mathcal{L}_{\mathcal{E}}$ of class (\ref{ed1}) which is finite-dimensional and is spanned by

\begin{equation}\label{eqalg} \begin{array}{l}\nonumber Y_1  =  \partial_t, \quad Y_2  =  \partial_x,\quad  Y_3  =  \partial_u,\\ \\

\nonumber Y_4  =  t \partial_t -f \partial_f -g \partial_g-\alpha \partial_{\alpha}-\beta \partial_{\beta}- \gamma \partial_{\gamma}-\phi \partial_{\phi}, \\ \\

\nonumber Y_5  =  x \partial_x +f \partial_f + 2 \alpha \partial_{\alpha}+3 \beta \partial_{\beta}+4 \gamma \partial_{\gamma}+2\phi \partial_{\phi}, \\ \\

\nonumber Y_6  =  u \partial_u +f \partial_f +g \partial_g+\phi \partial_{\phi}, \\ \\

\nonumber Y_7  =  t \partial_x +u \partial_f, \quad Y_8  =  \partial_{\alpha} -u \partial_{\phi}, \\ \\

\nonumber Y_9  =  \partial_f, \quad Y_{10}  =   \partial_{\phi}.

\end{array} \end{equation}

\noindent From generators (\ref{eqalg}) we get the finite form of the equivalence transformations:

\begin{theorem}\label{th1} The equivalence group of class (\ref{ed1}) consists of the transformations\\
\begin{equation}\label{equivgroup}
\begin{array}{l}
 \displaystyle \tilde{t}  =  \left( t+ \epsilon_1  \right) e^{\epsilon_4}, \quad \tilde{x}  =  \left(x+\epsilon_7 t+\epsilon_2 \right)e^{\epsilon_5}, \quad \tilde{u}  =  \left(u + \epsilon_3 \right) e^{\epsilon_6}, \\ \\
 \tilde{f}  =  \displaystyle \left( f+\epsilon_7 u+\epsilon_9 \right)e^{-\epsilon_4+\epsilon_5+\epsilon_6}, \quad \displaystyle \tilde{g}= g e^{-\epsilon_4+\epsilon_6},  \quad   \tilde{\alpha}  =  \displaystyle \left( \alpha +\epsilon_8 \right)e^{-\epsilon_4+2\epsilon_5}, \\ \\
 \tilde{\beta}   =   \displaystyle \beta e^{-\epsilon_4+3\epsilon_5},  \quad  \tilde{\gamma}  = \displaystyle \gamma e^{-\epsilon_4+4\epsilon_5}, \quad \tilde{\phi}  = \left( \phi -\epsilon_8 u+\epsilon_{10} \right)e^{-\epsilon_4+2\epsilon_5+\epsilon_6},
\end{array}
\end{equation}

\noindent where $\epsilon_i$, $i =1, \ldots, 10$, are arbitrary constants.\\

\end{theorem}

\noindent \textbf{Proof.}  From (\ref{equivgroup}) it follows that $u=e^{-\epsilon_6}\tilde{u}-\epsilon_3$. Thus,
\begin{equation}\label{derivatives}
u_t=\displaystyle  e^{\epsilon_4-\epsilon_6}\tilde{u}_{\tilde{t}}+ \epsilon_7 e^{\epsilon_5-\epsilon_6}\tilde{u}_{\tilde{x}}, \quad u_{x}=e^{\epsilon_5-\epsilon_6}\tilde{u}_{\tilde{x}}, \quad u_{xx}=e^{2\epsilon_5-\epsilon_6}\tilde{u}_{\tilde{x}\tilde{x}}, \quad u_{xxx}=e^{3\epsilon_5-\epsilon_6}\tilde{u}_{\tilde{x}\tilde{x}\tilde{x}}, \quad u_{xxxx}=e^{4\epsilon_5-\epsilon_6}\tilde{u}_{\tilde{x}\tilde{x}\tilde{x}\tilde{x}}.
\end{equation}

\noindent Moreover,

\begin{equation}\label{func1}
\begin{array}{l}
f=\displaystyle -\epsilon_7 u -\epsilon_9+ e^{\epsilon_4-\epsilon_5-\epsilon_6}\tilde{f}, \quad g= \displaystyle e^{\epsilon_4-\epsilon_6}\tilde{g}, \quad \alpha=  -\epsilon_8 + e^{\epsilon_4-2\epsilon_5}\tilde{\alpha}, \\ \\ \beta =e^{\epsilon_4-3\epsilon_5}\tilde{\beta}, \quad \gamma=e^{\epsilon_4-4\epsilon_5}\tilde{\gamma}, \quad \phi= \displaystyle \epsilon_8 u-\epsilon_{10}+e^{\epsilon_4-2\epsilon_5-\epsilon_6}\tilde{\phi}.
\end{array}
\end{equation}

\noindent Hence,

\begin{equation}\label{funcder}
f_u=\displaystyle  -\epsilon_7 + e^{\epsilon_4-\epsilon_5}\tilde{f}_{\tilde{u}}, \qquad \phi_u =\epsilon_8+e^{\epsilon_4-2\epsilon_5}\tilde{\phi}_{\tilde{u}}, \qquad \phi_{uu}=e^{\epsilon_4-2\epsilon_5+\epsilon_6}\tilde{\phi}_{\tilde{u}\tilde{u}}.
\end{equation}

\noindent Substituting (\ref{derivatives})-(\ref{funcder}) and then dividing the result by $e^{\epsilon_4-\epsilon_6}$ we obtain

$$ \tilde{u}_{\tilde{t}}+ \tilde{f}_{\tilde{u}} \tilde{u}_{\tilde{x}}+ \tilde{\alpha} \tilde{u}_{\tilde{x}\tilde{x}}+ \tilde{\phi}_{\tilde{u}\tilde{u}} \tilde{u}_{\tilde{x}}^2 + \tilde{\phi}_{\tilde{u}} \tilde{u}_{\tilde{x}\tilde{x}}+\tilde{\beta} \tilde{u}_{\tilde{x}\tilde{x}\tilde{x}}+\tilde{\gamma} \tilde{u}_{\tilde{x}\tilde{x}\tilde{x}\tilde{x}}=\tilde{g}. \quad {}_{\square}$$

\section{Theorem on Projections}

At this point, it is useful to know the most general symmetry algebra admitted by equation (\ref{ed1}) when the functions $f$, $g$, $\alpha$, $\beta$, $\gamma$ and $\phi$ are arbitrary. This symmetry algebra is known as the principal Lie algebra $\mathcal{L}_p$. Equivalence transformations are widely used in the study of PDEs, in particular, in nonlinear PDEs. For instance, they can be used for obtaining the principal Lie algebra as well as for finding extensions of the principal Lie algebra once the arbitrary elements are fixed. Let us introduce the following projections of the equivalence operator (\ref{generator})

\begin{equation}
\label{varproj} X= pr_{({\rm x}, u)}(Y)= \tau \partial_t + \xi \partial_x+\eta \partial_u,
\end{equation}

\begin{equation}
\label{funcproj} Z= pr_{(u, \Psi)}(Y)= \eta \partial_u+ \omega^1 \partial_f  + \omega^2 \partial_g + \omega^3 \partial_{\alpha} + \omega^4 \partial_{\beta} + \omega^5 \partial_{\gamma}+ \omega^6 \partial_{\phi},
\end{equation}

\noindent where ${\rm x}=(t,x)$ and $\Psi=(f,g,\alpha,\beta,\gamma,\phi)$ represent the set of independent variables and the set of arbitrary functions respectively. The equivalence generator $Y$ of class (\ref{ed1}) can be written as a linear combination of the elements of $\mathcal{L}_{\mathcal{E}}$

\begin{equation}
\begin{array}{rcl}
Y & = & \displaystyle \sum_{i=1}^{10} c_i Y_i \\ \\
& = & \left( c_1+ c_4 t \right) \partial_t + \left( c_2 + c_5 x+ c_7 t \right) \partial_x + \left( c_3 + c_6 u \right) \partial_u+ \left( c_9+ c_7 u+ \left(-c_4+c_5+c_6 \right) f \right) \partial_f \\ \\
& & +\left( -c_4+c_6 \right) g \partial_g+ \left( c_8 + \left( -c_4+2 c_5 \right) \alpha \right) \partial_{\alpha}+ \left( -c_4+3 c_5 \right) \beta \partial_{\beta} \\ \\

& & + \left( -c_4+4 c_5  \right) \gamma \partial_{\gamma} + \left(c_{10}-c_8 u+ \left( -c_4+2 c_5+c_6 \right) \phi \right) \partial_{\phi}.
\end{array}
\end{equation}
\medskip
\noindent In order to obtain the principal Lie algebra we apply the following theorem:

\begin{theorem}
\label{thproj1}
Let $Y$ be an equivalent operator for equation (\ref{ed1}), the projection $X= pr_{({\rm x}, u)}(Y) \in \mathcal{L}_{p}$ if and only if the projection $Z= pr_{(u, \Psi)}(Y)$ is identically zero.
\end{theorem}

\noindent Taking into account Theorem \ref{thproj1}, we require the vanishing of the projection $Z$

\begin{equation}
\begin{array}{rcl}
Z & = & \displaystyle  \left( c_3 + c_6 u \right) \partial_u+ \left( c_9+ c_7 u+ \left(-c_4+c_5+c_6 \right) f \right) \partial_f  +\left( -c_4+c_6 \right) g \partial_g+ \left( c_8 + \left( -c_4+2 c_5 \right) \alpha \right) \partial_{\alpha} \\ \\

& & + \left( -c_4+3 c_5 \right) \beta \partial_{\beta}+ \left( -c_4+4 c_5  \right) \gamma \partial_{\gamma} + \left(c_{10}-c_8 u+ \left( -c_4+2 c_5+c_6 \right) \phi \right) \partial_{\phi},
\end{array}
\end{equation}

\noindent this leads us to

$$ c_3= c_4= c_5= c_6= c_7= c_8= c_9 = c_{10}=0.$$

\noindent  Thus, we get that the principal Lie algebra $\mathcal{L}_{p}$ is given by

$$ X= pr_{({\rm x}, u)}(Y)= c_1 Y_1+ c_2 Y_2.$$

\noindent Consequently, for arbitrary functions $f$, $g$, $\alpha$, $\beta$, $\gamma$ and $\phi$, equation (\ref{ed1}) admits the two-dimensional Lie algebra generated by	
\begin{equation}\label{principalalg}
X_1= \partial_t, \qquad X_2= \partial_x.
\end{equation}

\noindent We are interested in obtaining extensions of the principal Lie algebra, for this purpose, we will make use of the following theorem:

\begin{theorem}
\label{thproj2}
Let $Y$ be an equivalent operator for equation (\ref{ed1}). The operator $X= pr_{({\rm x}, u)}(Y)$ is admitted by equation (\ref{ed1}) with specific functions
\begin{equation}\label{func} f=F(u), \quad g=G(u), \quad \alpha=A(u), \quad \beta= B(u), \quad \gamma= \Gamma(u), \quad \phi=\Phi(u), \end{equation}
if and only if these functions are invariant with respect to the projection $Z= pr_{(u, \Psi)}(Y)$.
\end{theorem}

\noindent The proof of Theorem \ref{thproj1} and \ref{thproj2} can be found in \cite{ibragimov2002}.\\

\section{Optimal system of one-dimensional subalgebras}

Generally, it is not usually feasible to obtain all the possible group-invariant solutions, since there can be an infinite number of Lie subgroups of the Lie group of symmetries $G$ of a given equation. Therefore, we would like to classify all the possible invariant solutions into different classes such as two solutions belonging to the same class are equivalent (one solution can be transformed into the other under the action of an element of the Lie symmetry group) and solutions belonging to different classes are not equivalent (these solutions are not related by any element of the Lie symmetry group). This classification problem can be solved by constructing the optimal system of subalgebras \cite{olver,ovsian}.\\

An optimal system of one-dimensional Lie subalgebras includes essential information about different types of invariant solutions \cite{olver,ovsian}. Let $\mathcal{G}$ be the Lie algebra of $G$. We say two one-dimensional subalgebras $h$ and $\tilde{h}$ are equivalent if $\tilde{h}=\mbox{Ad } g (h) $ where $\mbox{Ad } g$ is the adjoint action of $g$ on $\mathcal{G}$. Moreover, for each $V \in \mathcal{G}$ is defined a linear operator $\, \mbox{ad } V: \mathcal{G} \longrightarrow \mathcal{G}$, $\, \mbox{ad } V(W)=[V,W]$, where $[\,,\,]$ represents the Lie bracket. A set of one-parameter subalgebras forms an optimal system if every one-parameter subalgebra of $\mathcal{G}$ is equivalent to an exclusive element of the set under some component of the adjoint representation $\tilde{h}=\mbox{Ad } g (h), g \in G $.\\

Through the employment of the exponential map $\mbox{Exp}$ from $\mathcal{G}$ to $G$ it is possible to construct the adjoint action $\mbox{Ad } G$ of the underlying Lie group by means of the formula \cite{olver}

$$
\mbox{Ad}(exp(\epsilon V))W=\displaystyle \sum_{n=0}^{\infty} \frac{\epsilon^n}{n!}(\mbox{ad } V)^n(W)=W-\epsilon [V,W]+ \frac{\epsilon^2}{2} [V,[V,W]]-\ldots
$$

The classification of all nonequivalent equations (with respect to a given equivalence group ${\cal E}$) admitting extensions of the principal Lie algebra is known as a preliminary group classification. The method of preliminary group classification is simple and effective when the classification is based on finite-dimensional equivalence algebra $\mathcal{L}_{\mathcal{E}}$. In our case, we take into account the ten-dimensional subalgebra (\ref{eqalg}) and use it for obtaining a preliminary group classification. According to Theorem \ref{thproj2} the problem of preliminary group classification of equation (\ref{ed1}) with respect to the finite-dimensional algebra (\ref{eqalg}) is equivalent to constructing the optimal system of subalgebras given by

\begin{equation}\label{projections} \begin{array}{l}\nonumber Z_1  =  \partial_u, \quad Z_2  = -f \partial_f-g \partial_g -\alpha \partial_{\alpha}-\beta \partial_{\beta}- \gamma \partial_{\gamma}-\phi \partial_{\phi}, \quad Z_3  =  f \partial_f + 2 \alpha \partial_{\alpha}+3 \beta \partial_{\beta}+4 \gamma \partial_{\gamma}+2\phi \partial_{\phi}, \\ \\

\nonumber Z_4  =  u \partial_u +f \partial_f +g \partial_g+\phi \partial_{\phi}, \quad Z_5  = u \partial_f, \quad  Z_6  =  \partial_{\alpha} -u \partial_{\phi}, \quad Z_7  =  \partial_f, \quad Z_{8}  =   \partial_{\phi},

\end{array} \end{equation}

\noindent where $Z_i$, $i=1,\ldots,8$, are the nonzero projections of (\ref{eqalg}) on the space $(u,\Psi )$.\\

Therefore, we consider $\mathcal{G}$ the symmetry algebra with basis $\{ Z_1,\ldots,Z_8 \}$. First, we construct the commutator table for $\mathcal{G}$ which is a $8\times 8$ table
whose $(i,j)$-th entry expresses the Lie bracket $[Z_i,Z_j]$. Commutator table is shown in Table \ref{table1}. From the commutator table for this algebra, we construct the adjoint action $\mbox{Ad }G$ on the basis $\{ Z_1,\ldots,Z_8 \}$ of $\mathcal{G}$. We show the adjoint action in Table \ref{table2}, where the $(i,j)$-th entry gives $\mbox{Ad}(exp(\epsilon Z_i))Z_j$.\\

\begin{table}
\centering
\caption{Table of commutators}
\begin{tabular}{ccccccccc}
\hline
\hline
  & $Z_1$ & $Z_2$ & $Z_3$ & $Z_4$ & $Z_5$ & $Z_6$ & $Z_7$ & $Z_8$ \\
\hline
$Z_1$ & 0 & 0 & 0 & $Z_1$ & $Z_7$ & $-Z_8$ & 0 & 0 \\

$Z_2$ & 0 & 0 & 0 & 0 & $Z_5$ & $Z_6$ & $Z_7$ & $Z_8$ \\

$Z_3$ & 0 & 0 & 0 & 0 & $-Z_5$ & $-2 Z_6$ & $-Z_7$ & $-2 Z_8$ \\

$Z_4$ & $-Z_1$ & 0 & 0 & 0 & 0 & 0 & $-Z_7$ & $-Z_8$ \\

$Z_5$ & $-Z_7$ & $-Z_5$ & $Z_5$ & 0 & 0 & 0 & 0 & 0 \\

$Z_6$ & $Z_8$ & $-Z_6$ & $2 Z_6$  & 0 & 0 & 0 & 0 & 0 \\

$Z_7$ & 0 & $-Z_7$ & $Z_7$ & $Z_7$ & 0 & 0 & 0 & 0 \\

$Z_8$ & 0 & $-Z_8$ & $2 Z_8$ & $Z_8$ & 0 & 0 & 0 & 0 \\
\hline
\hline
\end{tabular}
\label{table1}
\end{table}

\begin{table}
\centering
\caption{The adjoint action on the basis $\{ Z_1,\ldots,Z_8 \}$}
\begin{tabular}{ccccccccc}
\hline
\hline
Ad  & $Z_1$ & $Z_2$ & $Z_3$ & $Z_4$ & $Z_5$ & $Z_6$ & $Z_7$ & $Z_8$ \\
\hline
$Z_1$ & $Z_1$ & $Z_2$ & $Z_3$ & $Z_4-\epsilon Z_1$ & $Z_5-\epsilon Z_7$ & $Z_6+\epsilon Z_8$ & $Z_7$ & $Z_8$ \\

$Z_2$ & $Z_1$ & $Z_2$ & $Z_3$ & $Z_4$ & $e^{-\epsilon}Z_5$ & $e^{-\epsilon}Z_6$ & $e^{-\epsilon}Z_7$ & $e^{-\epsilon}Z_8$ \\

$Z_3$ & $Z_1$ & $Z_2$ & $Z_3$ & $Z_4$ & $e^{\epsilon}Z_5$ & $e^{2\epsilon} Z_6$ & $e^{\epsilon}Z_7$ & $e^{2\epsilon} Z_8$ \\

$Z_4$ & $e^{\epsilon}Z_1$ & $Z_2$ & $Z_3$ & $Z_4$ & $Z_5$ & $Z_6$ & $e^{\epsilon}Z_7$ & $e^{\epsilon}Z_8$ \\

$Z_5$ & $Z_1+\epsilon Z_7$ & $Z_2+\epsilon Z_5$ & $Z_3-\epsilon Z_5$ & $Z_4$ & $Z_5$ & $Z_6$ & $Z_7$ & $Z_8$ \\

$Z_6$ & $Z_1-\epsilon Z_8$ & $Z_2+\epsilon Z_6$ & $Z_3-2\epsilon Z_6$  & $Z_4$ & $Z_5$ & $Z_6$ & $Z_7$ & $Z_8$ \\

$Z_7$ & $Z_1$ & $Z_2+\epsilon Z_7$ & $Z_3-\epsilon Z_7$ & $Z_4-\epsilon Z_7$ & $Z_5$ & $Z_6$ & $Z_7$ & $Z_8$ \\

$Z_8$ & $Z_1$ & $Z_2+\epsilon Z_8$ & $Z_3-2\epsilon Z_8$ & $Z_4-\epsilon Z_8$ & $Z_5$ & $Z_6$ & $Z_7$ & $Z_8$ \\
\hline
\hline
\end{tabular}
\label{table2}
\end{table}

In substance, the optimal system can be constructed by taking a general element $Z=a_1 Z_1+a_2 Z_2+\ldots+a_8 Z_8 \in \mathcal{G}$, where $a_i$, $i=1,\ldots,8$, are arbitrary constants, and simplifying it as much as possible through well-considered applications of adjoint maps to $Z$, $\mbox{Ad}(exp(\epsilon Z_i))Z$, and discrete symmetries. For further information on how optimal system can be constructed, one can refer, as example, to reference \cite{olver}.\\

\noindent To conclude, we obtain the following optimal system of one-dimensional subalgebras of (\ref{projections})

\begin{equation}\label{optimal system} \begin{array}{l} Z^{(1)}  = Z_1, \quad Z^{(2)}  = Z_2, \quad Z^{(3)}  =  Z_1+Z_2, \quad Z^{(4)}= Z_1+Z_8, \quad Z^{(5)}= a Z_2+Z_3, \quad Z^{(6)}=Z_7+a Z_8, \\ \\

 Z^{(7)}  =  Z_1+a Z_2+Z_3, \quad Z^{(8)}  = a Z_2+ b Z_3+Z_4, \quad  Z^{(9)} =  c Z_2+ d Z_3+Z_5, \quad Z^{(10)}  = a Z_2+b Z_3 +Z_6 , \\ \\

 Z^{(11)}= Z_2+Z_3+Z_7, \quad Z^{(12)}= 2 Z_2+ Z_3 +Z_8, \quad Z^{(13)}= Z_4+Z_5+ a Z_6, \quad Z^{(14)}= Z_4-Z_5+ a Z_6,\\ \\

 Z^{(15)}= Z_1+a Z_2+b Z_3+Z_5, \quad Z^{(16)}= -Z_1+a Z_2+b Z_3+Z_5, \quad Z^{(17)}=Z_1+a Z_2+b Z_3+Z_6, \\ \\

 Z^{(18)}= -Z_1+a Z_2+b Z_3+Z_6, \quad Z^{(19)}= Z_1+Z_2+Z_3+Z_7, \quad Z^{(20)}= Z_1+2 Z_2+Z_3+Z_8,\\ \\

Z^{(21)}= m Z_2+ m Z_3+Z_4+Z_5, \quad Z^{(22)}= 2 a Z_2+ a Z_3+ Z_4+ Z_6, \quad Z^{(23)}=(1+a)Z_2+a Z_3+Z_4+Z_7, \\ \\

Z^{(24)}=(1+2a) Z_2+ a Z_3+Z_4+Z_8, \quad Z^{(25)}=a Z_2+a Z_3+Z_6+Z_7, \quad Z^{(26)}=Z_2+Z_4+Z_7+m Z_8, \\ \\

Z^{(27)}=Z_1+2a Z_2+a Z_3+Z_5+m Z_6, \quad Z^{(28)}=-Z_1+2a Z_2+a Z_3+Z_5+m Z_6, \\ \\

Z^{(29)}=Z_1+a Z_2+a Z_3+Z_6+m Z_7, \quad Z^{(30)}=-Z_1+a Z_2+a Z_3+Z_6+m Z_7, \\ \\

Z^{(31)}=-Z_2- Z_3+ Z_4+Z_5+ Z_8, \quad Z^{(32)}=2 Z_2+ Z_3+ Z_4+m Z_6+ Z_7, \\ \\

Z^{(33)}=2 a Z_2+a Z_3+Z_5+b Z_6+Z_8, \quad  Z^{(34)}=2 a Z_2+a Z_3+Z_5+b Z_6-Z_8,\\ \\

\end{array} \end{equation}

\noindent where $a$, $b$, $c \neq 2d$, $m \neq 0$ are arbitrary constants.

\section{Reductions and exact solutions}

The symmetry group of a partial differential equation (PDE) is the largest transformation group
which acts on dependent and independent variables of the equation so it transforms solutions of
the equation into other solutions. Lie symmetry groups are considered to be one of the most powerful
methods to construct exact solutions of PDEs. Local symmetries admitted by a PDE are useful for obtaining invariant solutions. This technique is based on the following: if a differential equation is invariant under a Lie group of transformations, then a reduction transformation exists. In the case of PDEs with two independent variables, this technique yields a similarity variable and a similarity solution which allow us to transform the PDE into an ODE, which is generally easier to solve.\\

Lie point symmetries can be obtained by using the invariance criterion \cite{olver,ovsian} which leads to an overdetermined linear system of equations called determining system. The problem lies in the fact that, when arbitrary functions appear, the determining system usually becomes very difficult to solve. In these cases, the method of preliminary group classification is an effective way to give a group classification. This method does not guarantee a priori a complete group classification, however, it allows us to obtain symmetries when the determining system is intricate.\\

Thus, by applying Theorem \ref{thproj2} to the optimal system (\ref{optimal system}) it is possible to obtain all nonequivalent equations (\ref{ed1}) which admit extensions of the principal Lie algebra by one, namely, those equations belonging to (\ref{ed1}) such that they admit, besides the primary operators $X_1$ and $X_2$ (\ref{principalalg}), an additional operator $X_3$. We note that due to the extension of the optimal system, we do not show all nonequivalent equations. \\

From the associated equivalence algebra (\ref{eqalg}) we can separate equivalence generators for class (\ref{ed1}) into three different groups: translations $\left( Y_1, Y_2, Y_3, Y_9 \mbox{ and } Y_{10} \right)$, dilatations $\left( Y_4, Y_5 \mbox{ and } Y_{6} \right)$ and generalized Galilean transformations $\left( Y_7 \mbox{ and } Y_{8} \right)$. In order to illustrate the algorithm of passing from operators of the equivalence algebra (\ref{eqalg}) or, more specifically, from operators (\ref{optimal system}) to the form of the functions $f$, $g$, $\alpha$, $\beta$, $\gamma$ and $\phi$, and the additional operator $X_3$, we consider, by way of example, some generators involving representatives of the three different groups. Once we have obtained the new generator $X_3$, we determine the similarity variable and the similarity solution which allow us to transform the PDE which admits $X_3$ into an ODE. Finally, we obtain some exact solutions of the corresponding ODE.\\

\noindent \textbf{Example 1.} To begin with, we take the generator

\begin{equation}
\label{y3y4} Y= Y_3+Y_4=   t \partial_t+\partial_u -f \partial_f -g \partial_g-\alpha \partial_{\alpha}-\beta \partial_{\beta}- \gamma \partial_{\gamma}-\phi \partial_{\phi}. \\ \\
\end{equation}

\noindent Following Theorem \ref{thproj2}, we require the invariance of (\ref{func}) with respect to the operator

\begin{equation}
\label{zy3y4} Z^{(3)}= pr_{(u, \Psi)}(Y)= \partial_u -f \partial_f -g \partial_g-\alpha \partial_{\alpha}-\beta \partial_{\beta}- \gamma \partial_{\gamma}-\phi \partial_{\phi}.
\end{equation}

\noindent That means,

$$ \begin{array}{ccccc}
Z^{(3)}(f-F(u))\,{|}_{f=F(u)}=0, & \quad & Z^{(3)}(g-G(u))\,{|}_{g=G(u)}=0, & \quad &  Z^{(3)}(\alpha-A(u))\,{|}_{\alpha=A(u)}=0, \\ \\

Z^{(3)}(\beta-B(u))\,{|}_{\beta=B(u)}=0, & \quad & Z^{(3)}(\gamma-\Gamma(u))\,{|}_{\gamma=\Gamma(u)}=0, & \quad &  Z^{(3)}(\phi-\Phi(u))\,{|}_{\phi=\Phi(u)}=0,
\end{array}$$

\noindent which yields

$$ f= c_1 e^{-u}, \quad g=c_2 e^{-u}, \quad \alpha= c_3 e^{-u}, \quad \beta= c_4 e^{-u}, \quad \gamma= c_5 e^{-u}, \quad \phi=c_6 e^{-u},$$

\noindent where $c_i$ are arbitrary constants. Therefore, equation (\ref{ed1}) takes the following form

\begin{equation}
\label{eqy3y4} u_t+e^{-u} \left( -c_1 u_x+ \left( c_3 - c_6 \right) u_{xx}+c_6 u_x^{2}+c_4 u_{xxx}+c_5 u_{xxxx}- c_2 \right)=0,
\end{equation}

\noindent which admits the additional symmetry

\begin{equation}
\label{geny3y4} X_3= t \partial_t +\partial_u.
\end{equation}
The corresponding invariant solution is written in the form
\begin{equation}\label{transf1} z=x,\qquad u=\mbox{ln}\left|t \right| +h(z),\end{equation} where $h(z)$ must satisfy
\begin{equation}\label{ode1}{c_5}\,h''''+{c_4}\,h'''+{ \left( c_3 -c_6 \right)}\,h''+{c_6}\,\left(h'\right)^2-{ c_1}\,h'+e^{h}-{c_2}=0.\end{equation}

\noindent Since that $c_5\neq 0$, equation (\ref{ode1}) can be written as an autonomous equation of fourth order $h''''=F(h,h',h'',h''')$. The substitution $w(h)=(h')^2$ leads to the following third-order equation

\begin{equation}\label{subode1} 2 c_5 w w'''+ \left( c_5 w' \pm 2 c_4 \sqrt{w} \right) w''+2 \left(c_3-c_6\right) w'+4 c_6 w \mp 4 c_1 \sqrt{w}-4 c_2+4
   e^h=0. \end{equation}


\bigskip

\noindent \textbf{Example 2.} Let us consider generator
\begin{equation}
\label{y4y6} Y= Y_4+Y_6=   t \partial_t+u \partial_u -\alpha \partial_{\alpha}-\beta \partial_{\beta}- \gamma \partial_{\gamma}. \\ \\
\end{equation}

\noindent By applying Theorem \ref{thproj2}, we require the invariance of (\ref{func}) with respect to the projection

\begin{equation}
\label{zy4y6} Z= \displaystyle \left. Z^{(8)} \right|_{a=1,b=0}=  u \partial_u -\alpha \partial_{\alpha}-\beta \partial_{\beta}- \gamma \partial_{\gamma},
\end{equation}

\noindent which yields

$$ f= c_1, \quad g=c_2, \quad \alpha= \frac{c_3}{u}, \quad \beta= \frac{c_4}{u}, \quad \gamma= \frac{c_5}{u}, \quad \phi=c_6,$$

\noindent where $c_i$ are arbitrary constants. Thus, equation (\ref{ed1}) is given by

\begin{equation}
\label{eqy4y6} u u_t+c_3 u_{xx}+c_4 u_{xxx}+c_5 u_{xxxx}- c_2 u=0,
\end{equation}

\noindent which admits the additional symmetry

\begin{equation}
\label{geny4y6} X_3= t \partial_t +u\partial_u.
\end{equation}
The corresponding similarity variable and similarity solution are written in the form
\begin{equation}\label{transf2} z=x,\qquad u=t \, h(z),\end{equation} where $h(z)$ must satisfy
\begin{equation}\label{ode2}{c_5}\,h''''+{c_4}\,h'''+c_3 \,h''+h^2-{c_2} \,h=0.\end{equation}

\noindent In this case, equation (\ref{ode2}) does not admit any nontrivial Lie symmetry but the substitution $w(h)=(h')^2$ leads to the following third-order equation

\begin{equation}\label{subode2} 2 c_5 w w'''+ \left( c_5 w' \pm 2 c_4 \sqrt{w} \right) w''+2 c_3 w'+4 h^2-4 c_2 h=0. \end{equation}

\bigskip

\noindent \textbf{Example 3.} Now, we consider the generator

\begin{equation}
\label{y4y5y6} Y= Y_5+Y_6=   x \partial_x +u \partial_u +2 f \partial_f+ g \partial_g+2 \alpha \partial_{\alpha}+3 \beta \partial_{\beta}+4 \gamma \partial_{\gamma}+3 \phi \partial_{\phi}. \\ \\
\end{equation}

\noindent The projection (\ref{funcproj}) is given by

\begin{equation}
\label{zy5y6} Z= \displaystyle \left. Z^{(8)} \right|_{a=0,b=1}= u \partial_u +2 f \partial_f+ g \partial_g+2 \alpha \partial_{\alpha}+3 \beta \partial_{\beta}+4 \gamma \partial_{\gamma}+3 \phi \partial_{\phi}.
\end{equation}

\noindent Requiring the invariance of (\ref{func}) with respect to the operator (\ref{zy5y6}) we obtain that the functions are given by

$$ f= c_1 u^2, \quad g=c_2 u, \quad \alpha= c_3 u^2, \quad \beta= c_4 u^3, \quad \gamma= c_5 u^4, \quad \phi=c_6 u^3,$$

\noindent where $c_i$ are arbitrary constants. Thus, we conclude that equation (\ref{ed1}) with the above functions

\begin{equation}
\label{eqy5y6} u_t+2 u u_x \left( c_1+3 c_6 u_x \right) + u^2 u_{xx} \left( c_3 +3 c_6 \right) +c_4 u^3 u_{xxx}+c_5 u^4 u_{xxxx}- c_2 u=0,
\end{equation}

\noindent admits the additional symmetry

\begin{equation}
\label{geny5y6} X_3= x \partial_x +u \partial_u.
\end{equation}
The corresponding invariant solution is written in the form
\begin{equation}\label{transf3} z=t,\qquad u=x\,h(z),\end{equation} where $h(z)$ must satisfy
\begin{equation}\begin{array}{l}\label{ode3}h' =-6 c_6 h^3-2 {  c_1} h^2+{  c_2} h.\end{array}\end{equation}
If $c_6=0$, equation (\ref{ode3}) is the Bernoulli equation, thus we can obtain the corresponding solution $h$ of the ODE
$$h(z)=\frac{c_2 e^{c_2 z}}{2 c_1 e^{c_2 z}-e^{c_0 c_2}},$$
where $c_0$ is a constant of integration. In the case that $c_6 \neq 0$, equation (\ref{ode3}) is an Abel equation of the first kind, more specifically is a separable equation whose implicit solution is given by

$$\displaystyle \log \left(2 h \left(3 c_6 y+c_1\right)-c_2\right)-2 \log (h) +\frac{2 \,c_1 \arctan \left(\frac{6 c_6
   h+c_1}{\sqrt{-c_1^2-6 c_2 c_6}}\right)}{\sqrt{-c_1^2-6 c_2 c_6}}= -2 c_2 z+c_0,$$

\noindent where $c_0$ is an arbitrary constant. \\

\bigskip

\noindent \textbf{Example 4.} Finally, we take into account generator

\begin{equation}
\label{y3y7y8} Y= Y_3+Y_7+Y_8=   t  \partial_x + \partial_u +u \partial_f+ \partial_{\alpha}-u \partial_{\phi}, \\ \\
\end{equation}

\noindent whose projection is

\begin{equation}
\label{zy3y7y8} Z= \displaystyle \left. Z^{(27)} \right|_{a=0,m=1}=\partial_u +u \partial_f+ \partial_{\alpha}-u \partial_{\phi}.
\end{equation}

\noindent The invariance of (\ref{func}) with respect to (\ref{zy3y7y8}) leads us to

$$
f= \frac{u^2}{2}+c_1, \quad g=c_2, \quad \alpha= u+c_3, \quad \beta= c_4, \quad \gamma= c_5, \quad \phi=-\frac{u^2}{2}+c_6,
$$

\noindent where $c_i$ are arbitrary constants. Equation (\ref{ed1}) with the above specific functions

\begin{equation}
\label{eqy3y7y8} u_t+u u_x +c_3 u_{xx} -u_x^2 +c_4 u_{xxx}+c_5 u_{xxxx}- c_2=0,
\end{equation}

\noindent admits the following additional operator

\begin{equation}
\label{geny3y7y8} X_3= t \partial_x + \partial_u.
\end{equation}
The corresponding invariant solution is written in the form
\begin{equation}\label{transf4} z=t,\qquad u=\frac{x}{t}+h(z),\end{equation} where $h(z)$ must satisfy
\begin{equation}\begin{array}{l}\label{ode4}z^2h'-{c_2}\,z^2+h\,z-1=0.\end{array}\end{equation}
For this reduction, the ODE is a first order equation whose solution is
$$ h \left( z \right) = \displaystyle \frac{c_2 z}{2} +\frac{\ln  \left( z \right) +{c_0}}{z},$$ where $c_0$ is a constant of integration. Undoing transformation (\ref{transf4}), we obtain a solution of equation (\ref{eqy3y7y8}) which is given by

$$ u \left( t,x \right) = \displaystyle \frac{c_2 t}{2} +\frac{x+\ln  \left( t \right) +{c_0}}{t}.$$

\section{Conclusions}

In this paper we have considered a fourth-order nonlinear partial differential equation containing six arbitrary functions depending on the dependent variable. We have constructed the equivalence group of the class which is exploited to simplify the classifying equations. We have determined that class (\ref{ed1}) has a ten-dimensional algebra (\ref{eqalg}). We have used this algebra to perform a preliminary group classification. According to a theorem on projections, the problem of preliminary group classification of class (\ref{ed1}) is reduced to the construction of the optimal system of subalgebras of the nonzero projections of (\ref{eqalg}) on the space $(u,f,g,\alpha,\beta,\gamma,\phi )$. From the optimal system is possible to obtain all nonequivalent equations (\ref{ed1}) which admit an additional symmetry. In order to clarify the algorithm, we have taken some nonequivalent equations admitting extensions of the principal Lie algebra. Moreover, we have shown the corresponding coefficients $f,g,\alpha,\beta,\gamma,\phi$, and the extra symmetry operator. Finally, by using these symmetries, we have reduced the PDE into ODEs and some exact solutions have been derived.

\section*{Acknowledgments}

The authors gratefully acknowledge helpful discussions on equivalence transformations and the collaboration received by Dr. Rita Tracin\`a (Universit\`a degli studi di Catania). The authors also acknowledge the financial support from Junta de Andaluc\'ia group FQM-201, Universidad de C\'adiz and they express their sincere gratitude to the Plan Propio de Investigaci\'on de la Universidad de C\'adiz.


\end{document}